\documentclass[10pt]{article}

\usepackage{amsmath,amssymb,amsthm,amscd,latexsym}

\newtheorem{theorem}{Theorem}
\newtheorem{lemma}{Lemma}
\newtheorem{proposition}{Proposition}

\newtheorem{corollary}{Corollary}

\newtheorem{example}{Example}

\newtheorem{remark}{Remark}
\newenvironment{rem}
 {\begin{remark}\normalfont}
 {\end{remark}}

\evensidemargin=\oddsidemargin

\newcommand{\q}[2]{\genfrac{[}{]}{0pt}{}{#1}{#2}}

\newcommand{\di}{\#}

\begin{document}

\vspace{.3in}

\centerline{STRUCTURE OF $T$ MODULES AND RESTRICTED DUALS: THE CLASSICAL}

\vspace{.3in}

\centerline{AND THE QUANTUM CASE}

\vspace{.4in}

\renewcommand{\thefootnote}{\fnsymbol{footnote}}

\centerline{Dijana Jakeli\'c\footnote{Current address: Department of Mathematics,
University of California, Riverside, CA 92521. E-mail: \texttt{jakelic@math.ucr.edu}}}

\renewcommand{\thefootnote}{\arabic{footnote}}

\vspace{.3in}

\centerline{Department of Mathematics, University of Virginia, Charlottesville, VA 22904}

\vspace{.6in}


\centerline{ABSTRACT}

\centerline{A concrete realization of Enright's $T$ modules is obtained. This is
used to show}

\centerline{their self-duality. As a consequence, the restricted duals of Verma modules
are also identified.}

\vspace{.6in}


\centerline{1 INTRODUCTION}

\vspace{.3in}

In \cite{E}, T. Enright introduced a completion functor on the category $\mathcal{I}(\mathfrak{g})$
of representations of a finite dimensional complex semisimple Lie algebra $\mathfrak{g}$.
Using this functor, he constructed algebraically the fundamental series representations.
In \cite{D}, V. Deodhar realized the functor via Ore localization and thus gave an explicit way
of constructing completions which also enabled him to prove Enright's uniqueness conjecture arising
in successive completions. Later in \cite{Jo}, A. Joseph generalized the functor to the
Bernstein - Gelfand - Gelfand category $\mathcal{O}(\mathfrak{g})$ and gave a refinement of the
Jantzen conjecture.

As completion is associated with a simple root, the whole process is built up from the $sl_2$-case.
In particular, Enright introduced $T$ modules for $sl_2$ and showed that, together
with Verma modules, they comprise all indecomposable objects of the category
$\mathcal{I}(sl_2)$ and every module in this category is a direct sum (not necessary finite)
of indecomposable ones. More recently, Y.~M. Zou \cite{Z} extended the notion of completion,
Deodhar's approach, and results related to $T$ modules to the quantum case when $q$
is not a root of unity. These extensions go through smoothly with necessary adjustments
to accomodate $q$-integers.

In this paper, we are concerned with an explicit structure of $T$ modules.
The necessity for such description arose
in our study \cite{J} of the interaction between completions and crystal bases introduced by
M. Kashiwara \cite{Ka}.
A concrete realization of $T$ modules is obtained in a self-contained manner
by constructing their bases via PBW basis of $U_q(sl_2)$ (Theorem \ref{t:T}).
Consequently, the structure of any module in
$\mathcal{I}\left(U_q(sl_2)\right)$ becomes completely transparent.
Moreover, Theorem \ref{t:T} provides an elementary proof of
[1, Proposition 3.10] and [4, Proposition 4.3].

We study the restricted duals of all objects in the category
$\mathcal{I}\left(U_q(sl_2)\right)$ for which this notion is defined (Section 3).
In particular, the self-duality of $T$ modules is naturally revealed from
the explicit structure obtained in Theorem \ref{t:T}. The dual structure of Verma
modules then follows.

Although we present the aforementioned results for the quantum case, the proofs can be carried out similarly for the classical case. The corresponding
statement of Theorem \ref{t:T} in the classical case can be obtained, as usual, by setting $k=\mathbb{C}(q)$
with $q$ transcendental over $\mathbb{C}$, $t=q^h$, $\epsilon=1$, and letting $q \to 1$.
The corresponding statements of the remaining results are equally self-evident.

\vspace{.5in}


\centerline{2 CATEGORY $\mathcal{I}$ AND $T$ MODULES}

\vspace{.3in}

Let $\mathfrak{g}$ be a finite dimensional complex semisimple Lie algebra and $(a_{ij})_{1 \le i,j \le l}$ the corresponding
Cartan matrix.
There exist unique positive integers $d_i, \ i=1, \dots, l$, such that $\mbox{gcd} \ (d_1, \dots, d_l)=1$ and
the matrix $(d_ia_{ij})_{1 \le i,j \le l}$ is symmetric. Let $k$ be a field and $q \in k$ such that $q \ne 0$ and
$q^{2d_i} \ne 1 \ (1 \le i \le l)$. Denote $q_i=q^{d_i}$. The quantized enveloping algebra $U_q(\mathfrak{g})$ is the algebra
over $k$ with generators $e_i, f_i, t_i, t_i^{-1}, \ 1 \le i \le l$, and defining relations
\begin{align*}
 &t_it_i^{-1}=1=t_i^{-1}t_i, \ t_it_j=t_jt_i,\\
 &t_ie_jt_i^{-1}=q_i^{a_{ij}}e_j,  \ t_if_jt_i^{-1}=q_i^{-a_{ij}}f_j,\\
 &e_if_j-f_je_i=\delta_{ij} \frac{t_i-t_i^{-1}}{q_i-q_i^{-1}},\\
 &\sum_{s=0}^{1-a_{ij}} (-1)^s e_i^{(s)}e_je_i^{(1-a_{ij}-s)}=0, \ \sum_{s=0}^{1-a_{ij}} (-1)^s f_i^{(s)}f_jf_i^{(1-a_{ij}-s)}=0 \ \ (i \ne j)
\end{align*}
where as usual $\displaystyle{e_i^{(n)}=\frac{e_i^n}{[n]_{i}!}, \ f_i^{(n)}=\frac{f_i^n}{[n]_{i}!}}$,
$[n]_i ! = [1]_i [2]_i \dots [n]_i \ (n \in \mathbb{N})$, and
$\displaystyle{[n]_i = \frac{q_i^n-q_i^{-n}}{q_i-q_i^{-1}} \ (n \in \mathbb{Z})}$.

In particular, the following commutation relations hold for $s \in \mathbb{N}$:
\begin{equation}\label{commef}
 e_if_i^s=f_i^se_i+[s]_if_i^{s-1}[t_i;1-s], \ \  f_ie_i^s=e_i^sf_i-[s]_ie_i^{s-1}[t_i;s-1]
\end{equation}
where $\displaystyle{[t_i;a] = \frac{q_i^at_i-q_i^{-a}t_i^{-1}}{q_i-q_i^{-1}} \ \ \mbox{for} \ a \in \mathbb{Z}}$.

As an immediate consequence of the definition of $q$-integers, we have:

\begin{lemma}\label{l:[]}
 $[a+k]_i[b+k]_i-[a]_i[b]_i=[k]_i[a+b+k]_i \ \mbox{for all integers} \, a,b, \, \mbox{and} \ k$.
\end{lemma}

Denote by $U_q^+(\mathfrak{g})$ (resp. $U_q^-(\mathfrak{g})$) the subalgebra of $U_q(\mathfrak{g})$
generated by $e_i$ (resp. $f_i$), $1 \le i \le l$.

Henceforth, assume $\mbox{char} \, k \ne 2,3$ and $q$ is not a root of unity in $k$.
We recall few definitions following \cite{E} and \cite{Z}.

Let $\mathcal{I}\left( U_q(\mathfrak{g}) \right)$ be the category of
$U_q(\mathfrak{g})$-modules $M$
satisfying
(i) $M$ is a weight module,
(ii) $U_q^-(\mathfrak{g})$-action on $M$ is torsion free, and
(iii) $M$ is $U_q^+(\mathfrak{g})$-finite, i.e. $e_i$ acts locally nilpotently
on $M$ for all $i$.

Fix $i \in \{1, \dots, l\}$. Set
$M_\lambda = \{m \in M| \ t_im=\lambda m\} \ \mbox{for} \ \lambda \in k^\times = k \smallsetminus
\{0\}, \
M^{e_i} = \{m \in M| \ e_im=0\} , \ \mbox{and} \
M_\lambda^{e_i} = M_\lambda \cap M^{e_i}$.
A module $M$ in $\mathcal{I}\left( U_q(\mathfrak{g}) \right)$ is
said to be complete with respect
to $i$ if $ f_i^{n+1}: \ M_{\epsilon q_i^n}^{e_i} \to  M_{\epsilon q_i^{-n-2}}^{e_i}$
is bijective for all $n \in \mathbb{N}_0$ and $\epsilon = \pm 1$.
A module $N$ in $\mathcal{I}\left( U_q(\mathfrak{g}) \right)$ is
a completion of
$M$ with respect to $i$ provided
(i) $N$ is complete with respect to $i$,
(ii) $M$ is imbedded in $N$, and
(iii) $N \slash M$ is $f_i$-finite.

Now, we can consider $\mathfrak{g}=sl_2$. For brevity, we write $U_q=U_q(sl_2)$ and
$\mathcal{I}=\mathcal{I}\left(U_q(sl_2)\right)$. Henceforth, $\epsilon = \pm 1$ and
$\displaystyle{[n]=\frac{q^n-q^{-n}}{q-q^{-1}}}$, etc.
Note that the Verma modules $M(\lambda)$ with $\lambda \in k^\times \smallsetminus
\{\epsilon q^{-n-2} | \, n \in \mathbb{N}_0\}$ are complete and the completion of
$M(\epsilon q^{-n-2}) \, (n \in \mathbb{N}_0)$ is $M(\epsilon q^n)$.
The quantum Casimir element
$\displaystyle{C = \frac{qt+q^{-1}t^{-1}}{(q-q^{-1})^2}} + fe$
acts on each $M(\lambda) \ (\lambda \in k^\times)$
as multiplication by the scalar $\displaystyle{\frac{q\lambda+q^{-1}\lambda^{-1}}{(q-q^{-1})^2}}$.
For $n \in \mathbb{N}_0$,
the left ideals $I(n,\epsilon)$ and $J(n,\epsilon)$ of $U_q$
are defined as follows:
\[
 I(n,\epsilon)=U_q \{t-\epsilon q^{-n-2}, \, e^{n+2}, \, (C-\epsilon
c)^2\} \ \mbox{and} \ J(n,\epsilon)=U_q \{t-\epsilon q^{-n-2}, \, e^{n+2}\}
\]
where $\displaystyle{c=\frac{q^{n+1}+q^{-n-1}}{(q-q^{-1})^2}}$.
We consider the $U_q$-modules
$T(n,\epsilon)=U_q \slash {I(n,\epsilon)}$ and
$S(n,\epsilon)=
U_q\slash{J(n,\epsilon)}$, and call $T(n,\epsilon)$ a $T$ module.


Similarly as in the classical case (cf. \cite{E}):

\begin{lemma}\label{l:A}
Let $L$ be the $k$-span of $\{ f^i e^j | \  j \le n+1, \ i,j
\in \mathbb{N}_0 \}$.
Then $U_q = L \oplus J(n,\epsilon)$ and $S(n,\epsilon)$ is a free
$U_q^{-}$-module isomorphic to $L$.
\end{lemma}

An explicit description of the structure of $T$ modules
is given with the following theorem.

\begin{theorem}\label{t:T}
The $U_q$-module
$T(n,\epsilon)$ is generated by $\displaystyle{z = \sum_{i=0}^{n} \epsilon^{n-i}
\frac{[n]! [n-i]!}{[i]!} f^i e^i}$.
Set $v_{-1}=0$, $z_n=0$ and, for $i \ge 0$,
$v_i = f^i e^{n+1}$ and $z_{n+1+i} = f^i z$. Then
$\{v_i, \, z_{n+1+i} | \, i \ge 0\}$ is a $k$-basis of $T(n,\epsilon)$. Moreover,
\begin{align*}
 t z_i &= \epsilon q^{n-2i} z_i\\
 f z_i &= z_{i+1}\\
 e z_i &= \epsilon [i][n-i+1] z_{i-1} + v_{i-1}
\end{align*}
for $i \ge n+1$
and $\{v_i | \, i \ge 0\} $ is a $k$-basis of the Verma submodule of $T(n,\epsilon)$
with highest weight $\epsilon q^n$
and the usual $U_q$-action.
\end{theorem}

\begin{proof}
Denote the image of $u \in L \subseteq U_q$ in $S(n,\epsilon)=U_q/J(n,\epsilon)$
also by $u$ (see Lemma \ref{l:A}).
We consider the action of $t$ and $C$ on $S(n,\epsilon)$. Since $t \cdot 1 =
\epsilon q^{-n-2} 1$, then for $i,j \in \mathbb{N}_0, \ j \le n+1$,
\begin{equation}\label{T1}
 t \cdot f^ie^j = tf^ie^j \cdot 1 = \epsilon q^{2(j-i)-n-2}f^ie^j.
\end{equation}
Clearly $S(n,\epsilon) = \oplus_{i \ge 0} S(n,\epsilon)_{\epsilon
q^{n-2i}}$ and $\dim S(n,\epsilon)_{\epsilon q^{n-2i}} \le n+2$ for all $i$.
For $0 \le j \le n+1$, it follows that $C \cdot e^j = \epsilon c_je^j + fe^{j+1}$ where
$\displaystyle{
c_j = \frac{q^{2j-n-1}+q^{-2j+n+1}}{(q-q^{-1})^2}}$
and, since $C$ is a central element,
$C \cdot f^ie^j = \epsilon c_jf^ie^j+f^{i+1}e^{j+1}$
for $i \in \mathbb{N}_0$.
Evidently $c_{n+1-j}=c_j$ for all $j$. Moreover
$c-c_j=[j][n-j+1]$ and therefore $c_j=c$ iff $j=0$ or $j=n+1$.

Now, we consider $S(n,\epsilon)^{\epsilon c}$, the submodule of $S(n,\epsilon)$
where $C-\epsilon c$ is locally nilpotent.
Since
$S(n,\epsilon)_{\epsilon q^{n-2i}}$ is invariant under $C$
for every $i \ge 0$, it suffices to look for the
generalized $\epsilon c$-eigenspace $\left(S(n,\epsilon)_{\epsilon q^{n-2i}}\right)_{(\epsilon c)}$
of $C$ in each
$S(n,\epsilon)_{\epsilon q^{n-2i}}$.
For $0 \le i \le n$, $C$ is given in a $k$-basis $\{e^{n-i+1},fe^{n-i+2}, \dots,
f^ie^{n+1}\}$ of $S(n,\epsilon)_{\epsilon q^{n-2i}}$ by

\begin{equation}\label{T3}
\begin{pmatrix}
  \epsilon c_{n-i+1} &0                     &\dots     &               &0\\
  1                  &\epsilon c_{n-i+2}    &\ddots    &               &\vdots\\
  0                  &1                     &\ddots    &               &\\
  \vdots             &\ddots                &\ddots    &\epsilon c_n   &0\\
  0                  &\dots                 &0         &1              &\epsilon c_{n+1}
\end{pmatrix}
\end{equation}
Its characteristic polynomial is
$\mbox{char}_{C}(t)=(t-\epsilon c)p_1(t)$ with $p_1(\epsilon
c) \ne 0$, and so
$\left(S(n,\epsilon)_{\epsilon q^{n-2i}}\right)_{(\epsilon c)}$
is $k$-spanned by $f^ie^{n+1}$, an $\epsilon c$-eigenvector of $C$.
Similarly, $S(n,\epsilon)_{\epsilon q^{n-2i}}$ for $i \ge n+1$ has a
$k$-basis $\{f^{i-n-1}, f^{i-n}e, \dots, f^ie^{n+1}\}=
f^{i-n-1} \cdot (\mbox{basis of} \ S(n,\epsilon)_{\epsilon q^{-n-2}})$
and $\mbox{char}_{C}(t)=(t-\epsilon c)^2p_2(t)$ with
$p_2(\epsilon c) \ne 0$.
Since $C$ is central,
$\left(S(n,\epsilon)_{\epsilon q^{n-2i}}\right)_{(\epsilon c)}$
has a $k$-basis $\{f^ie^{n+1}, f^{i-n-1}z\}$ for some
generalized eigenvector $z \in S(n,\epsilon)_{\epsilon q^{-n-2}}$.
Note that $z=\sum_{i=0}^n \alpha_if^ie^i$ with
$\displaystyle{\alpha_i = \epsilon^{n-i} \frac{[n]![n-i]!}{[i]!}}$
does the job. Namely,
$(C-\epsilon c)z
=\sum_{i=0}^n \alpha_i \left(\epsilon(c_i-c)f^ie^i+f^{i+1}e^{i+1}\right)
=\sum_{i=1}^n \left(-\alpha_i\epsilon[i][n-i+1]+\alpha_{i-1}\right)f^ie^i
+f^{n+1}e^{n+1}=f^{n+1}e^{n+1}$ and hence $(C-\epsilon c)^2z=0$.
Therefore, a $k$-basis of
$S(n,\epsilon)^{\epsilon c}$ is
\begin{equation}\label{T7}
 \{f^ie^{n+1}, \, f^iz | \, i \ge 0\}.
\end{equation}
Consequently,
\begin{equation}\label{Cc}
(C-\epsilon c)^2 \equiv 0 \ \mbox{on} \ S(n,\epsilon)^{\epsilon c}.
\end{equation}

The action of $e$ on the basis vectors from (\ref{T7}) follows from
(\ref{T1}) and commutation relations (\ref{commef}). Namely,
$ef^{i}e^j=f^{i}e^{j+1}-\epsilon[i][n+1+i-2j]f^{i-1}e^j$ for $i \ge 0$,
$0 \le j \le n+1$
(here: $f^{-1}=0$) and in addition, by
Lemma \ref{l:[]} and the definition of $\alpha_i$,
$ef^jz = - \epsilon [j][n+1+j]f^{j-1}z+f^{n+j}e^{n+1}$ for $j \ge 0$.
In particular, $ef^ie^{n+1}=\epsilon[i][n-i+1]f^{i-1}e^{n+1}$ for $i \ge 0$ and
$ez=f^ne^{n+1}$.

Denote the image of $1 \in U_q$ in $T(n,\epsilon)=U_q/I(n,\epsilon)$
also by $1$. Similarly as in \cite{E}, observe $J(n,\epsilon) \subseteq I(n,\epsilon)$
and let $\varphi : S(n,\epsilon) \to T(n,\epsilon)$ be the induced $U_q$-module
homomorphism with $\varphi(1)=1$. Note that $\varphi$ induces a surjection
$\Phi : S(n,\epsilon)^{\epsilon c} \to T(n,\epsilon)$ because $C-\epsilon c$ is
locally nilpotent on $T(n,\epsilon)$. Due to (\ref{Cc}), $\Phi$ is injective,
thus a $U_q$-module isomorphism. Denoting the image of
$u \in S(n,\epsilon)^{\epsilon c}$ under $\Phi$ again by $u$ should cause no
confusion. Setting $v_i=f^ie^{n+1}$ and $z_{n+1+i}=f^iz$ in $T(n,\epsilon)$
for $i \ge 0$ and using the previous calculations, we obtain in particular
$ez_i=\epsilon[n-i+1][i]z_{i-1}+v_{i-1}$ for $i \ge n+1$.
The remaining claims of the theorem are now evident.
\end{proof}

Schematically, $T(n,\epsilon)$ looks like
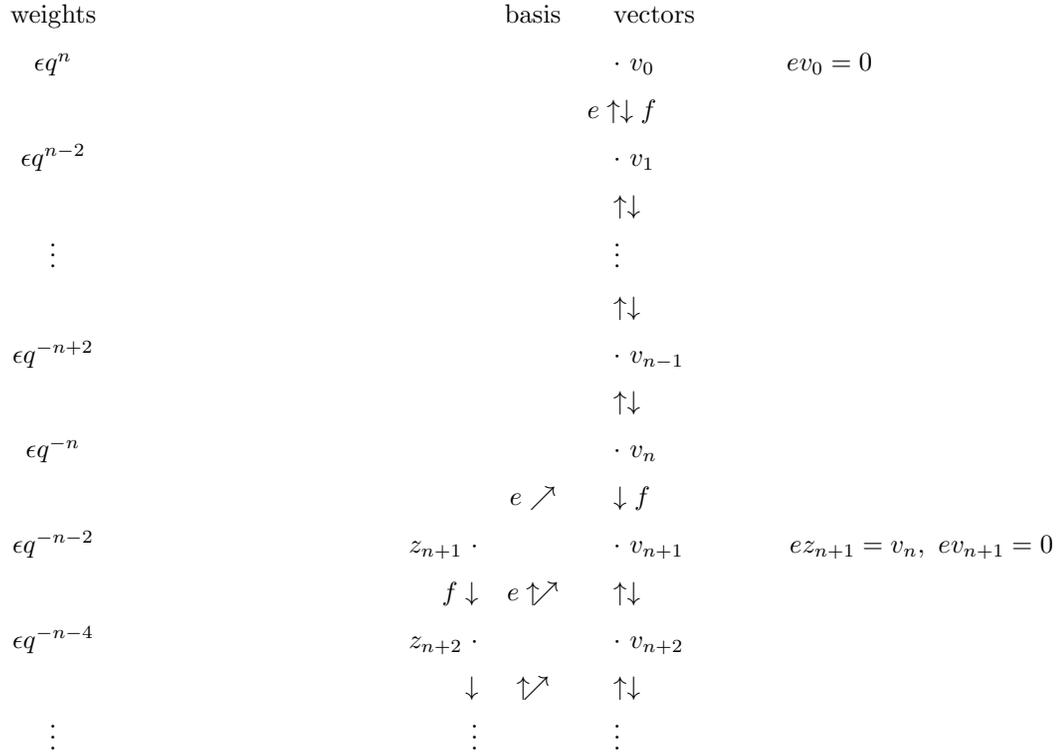
\begin{equation*}
 \begin{array}{crcl}
   \mbox{weights}           &\hspace{1.5in}
&\mbox{basis}                     &\quad \mbox{vectors}\\
   \epsilon  q^n            &\hspace{1in}                         &
                     &\quad \cdot \ v_0 \qquad \ \ \ \qquad ev_0=0\\
                            &                                     &
                     &e \uparrow \downarrow f\\
   \epsilon q^{n-2}         &\hspace{1.5in}                       &
                     &\quad \cdot \ v_1\\
                            &                                     &
                     &\quad \uparrow \downarrow\\
   \vdots                   &\hspace{1.5in}                       &
                     &\quad \vdots\\
                            &                                     &
                     &\quad \uparrow \downarrow\\
   \epsilon q^{-n+2}        &\hspace{1.5in}                       &
                     &\quad \cdot \ v_{n-1}\\
                            &                                     &
                     &\quad \uparrow \downarrow\\
   \epsilon q^{-n}          &\hspace{1.5in}                       &
                     &\quad \cdot \ v_n\\
                            &                                     &e \nearrow
                     &\quad \downarrow f\\
   \epsilon q^{-n-2}        &\hspace{1.5in}  z_{n+1} \ \cdot      &
                     &\quad \cdot \ v_{n+1} \qquad \qquad ez_{n+1}=v_n, \ ev_{n+1}=0\\
                            &f \downarrow                         &e
\uparrow\negthickspace\nearrow &\quad \uparrow \downarrow\\
   \epsilon q^{-n-4}        &\hspace{1.5in}  z_{n+2} \ \cdot      &
                     &\quad \cdot \ v_{n+2}\\
                            &\downarrow
&\uparrow\negthickspace\nearrow   &\quad \uparrow \downarrow\\
   \vdots                   &\hspace{1.5in}  \vdots               &
                     &\quad \vdots
 \end{array}
\end{equation*}

The following is immediate.

\begin{corollary} \label{MTM}
([1, Proposition 3.10], [4, Proposition 4.3]) \
The $U_q$-modules $S(n,\epsilon)$ and $T(n,\epsilon)$ belong to the
category $\mathcal{I} = \mathcal{I}\left(U_q(sl_2)\right)$,
$T(n,\epsilon) \cong S(n,\epsilon)^{\epsilon c}$,
$T(n,\epsilon)$ is complete and indecomposable, and there exists an
exact sequence
\[
 0 \to M(\epsilon q^n) \to T(n,\epsilon) \to M(\epsilon q^{-n-2}) \to 0.
\]
Moreover, if $M$ is a $U_q$-module in $\mathcal{I}$ such that
$M=M^{\epsilon c}$ and $v \in M_{\epsilon q^{-n-2}}$, then the map
$x \mapsto xv$ of $U_q$ into $M$ factors through $T(n,\epsilon)$.
\end{corollary}

Now, notice that
$v_{n+1} \in T(n,\epsilon)$ generates the Verma module
$M(\epsilon q^{-n-2})$,
and consider the $U_q$-module
$T(n,\epsilon) \slash M(\epsilon q^{-n-2})$.

\begin{proposition} \label{VTM}
There is a short exact sequence
\[
 0 \to V(n,\epsilon) \to T(n,\epsilon) \slash M(\epsilon q^{-n-2}) \to
M(\epsilon q^{-n-2}) \to 0
\]
where $V(n,\epsilon)$ is the irreducible $U_q$-module of dimension $n+1$ and
highest weight $\epsilon q^n$.
\end{proposition}
\begin{proof}
This is immediate from the filtration
$T(n,\epsilon) \supset M(\epsilon q^n) \supset M(\epsilon q^{-n-2}) \supset 0$,
the point being that
$T(n,\epsilon) \slash M(\epsilon q^n) \cong M(\epsilon q^{-n-2})$
and
$M(\epsilon q^n) \slash M(\epsilon q^{-n-2}) \cong V(n,\epsilon)$.
\end{proof}

The following result will be used only
at the end of Section 3 to keep the paper as self-contained
as possible.

\begin{proposition} \label{FG}([1, Proposition 3.11], [4, Proposition 4.5])
(i) \ The $M(\lambda) \ (\lambda \in k^\times)$ and the
$T(n,\epsilon) \ (n \in \mathbb{N}_0)$
are precisely all the indecomposable objects of the
category $\mathcal{I}$.

(ii) \ Every module in $\mathcal{I}$ is a direct sum (not necessarily finite) of indecomposable
ones.
\end{proposition}

\begin{rem}
In \cite{Z}, $k=\mathbb{C}(q)$ was used for convenience. No restriction on the field $k$ is needed other than
$\textnormal{char} \, k \ne 2$.
\end{rem}

\vspace{.5in}


\centerline{3 CATEGORY $\mathcal{O}$ AND RESTRICTED DUALS}

\vspace{.3in}

Let $\mathcal{O}\left(U_q(\mathfrak{g})\right)$ (cf. \cite{BGG,K,Z})
denote the category
consisting of $U_q(\mathfrak{g})$-modules $M$ such that
(1) $M$ is a weight module, i.e. $M=\oplus_{\omega \in \Omega}M_{\omega}$ for some $\Omega \subset
(k^\times)^l$,
(2) $\dim M_{\omega} < \infty$, and
(3) $P(M) \subset \cup_{1 \le i \le s} D(x_i)$ for some $s \in \mathbb{N}$
and $x_i \in (k^\times)^l \ (1 \le i \le s)$ where
$D(x_i)=\{y \in (k^\times)^l | \ y \le x_i\}$, i.e.
weights are contained in finitely many cones.

Evidently $U_q(\mathfrak{g})$ has an involutory antiautomorphism $\sigma$ such that
$\sigma (e_i)=f_i, \ \sigma (f_i)=e_i$ and $\sigma (t_i)=t_i$.
The following is the $q$-analogue of [7, Proposition 4.6].
Define the restricted dual $M^{res}$ of $M$ in
$\mathcal{O}\left(U_q(\mathfrak{g})\right)$ by
$M^{res} = \{f \in M^* | \ f(M_{\omega})=0 \ \mbox{for all but finitely many
weights} \ \omega \ \mbox{of} \ M\}$
where $M^* = \textnormal{Hom} \ (M,k)$.
$U_q(\mathfrak{g})$ acts on $M^{res}$ via $\sigma$:
\[
 (u f)(m)=f\left(\sigma (u) m\right) \, \mbox{for all} \, m \in M, u \in
U_q(\mathfrak{g}), f \in M^{res}.
\]
\noindent Denote $M^{res}$ with this action by $M^{\sigma}$.

As usual the formal character of
$M$
is defined as $\mbox{ch} \ M=\sum_{\omega \in \Omega}(\dim
M_\omega)e^\omega$ where $e^\omega$ is defined formally
to be the basis element of the group algebra $\mathbb{Z}[\Omega]$
corresponding to $\omega \in \Omega$.

\begin{proposition} \label{p:DGK}
Let $M \in Ob \ \mathcal{O}\left(U_q(\mathfrak{g})\right)$. Then
\begin{align*}
 &(i) \ \ \ M^\sigma \, \mbox{admits a weight space decomposition} \,
M^\sigma=\oplus_{\omega \in \Omega} M_{\omega}^\sigma \, \mbox{with} \,
M_{\omega}^\sigma=\{f \in M^* | \ f(M_{\mu})=0 \, \mbox{if} \, \mu
\ne \omega \}\\
 &(ii) \ \ M_{\omega}^\sigma \cong (M_{\omega})^* \, \mbox{as} \ k-\mbox{vector
spaces}\\
 &(iii) \ M^\sigma \in Ob \ \mathcal{O}\left(U_q(\mathfrak{g})\right)\\
 &(iv) \ \ \textnormal{ch} \ M^\sigma = \textnormal{ch} \ M\\
 &(v) \ \ \ M \cong (M^\sigma)^\sigma\\
 &(vi) \ \ M \to M^\sigma \, \mbox{is an exact contravariant functor of} \
\mathcal{O}\left(U_q(\mathfrak{g})\right) \, \mbox{onto itself}\\
 &(vii) \ V(\omega)^\sigma \cong V(\omega) \, \mbox{for all} \ \omega \in \Omega
\, \mbox{where} \, V(\omega) \, \mbox{is the irreducible} \
U_q(\mathfrak{g})-module \ \mbox{with highest weight} \ \omega\\
 &(viii) \mbox{If} \, M=\oplus_{i \in I} M_i \, \mbox{where} \,
 M_i \in Ob \ \mathcal{O}\left(U_q(\mathfrak{g})\right), \, \mbox{then} \,
 M^\sigma \cong \oplus_{i \in I} M_i^\sigma.
\end{align*}
\end{proposition}

\begin{proof}
The statements $(i)-(vii)$ have basically the same proofs
as in the classical case (cf. [7, Proposition 4.6], [10, Proposition 2.6.16]).

\noindent $(viii)$
Define $N_i=\{f \in M^\sigma | \, f(M_j)=0 \, \mbox{for} \, j \ne i\}$.
Then $M^\sigma = \oplus_{i \in I} N_i \cong \oplus_{i \in I} M_i^\sigma$.
\end{proof}

Next, we find the restricted duals of Verma modules and $T$ modules in
$\mathcal{O}=\mathcal{O}\left(U_q(sl_2)\right)$.


Let $\lambda \in k^\times$. The Verma module $M(\lambda)$
has a basis $\{v_i\}_{i \ge 0}$ where
$v_i=f^iv_0$. Then
$ tv_i=q^{-2i}\lambda v_i$,
$fv_i=v_{i+1}$, and
$\displaystyle{ev_i=[i] \ \frac{q^{1-i}\lambda-q^{i-1}{\lambda}^{-1}}{q-q^{-1}} \ v_{i-1} \
 (v_{-1}=0)}$.
By Proposition \ref{p:DGK}, \, $M(\lambda)^{\sigma}=\oplus_{\mu}
M(\lambda)_{\mu}^{\sigma}=\oplus_{i \ge 0}
M(\lambda)_{q^{-2i}\lambda}^{\sigma}$ where
$M(\lambda)_{q^{-2i}\lambda}^{\sigma}=\{f \in M(\lambda)^* | \
f(\left(M(\lambda)_{\mu}\right)=0$ for all $\mu \ne
 q^{-2i}\lambda\}$.
Since $M(\lambda)_{q^{-2i}\lambda}^{\sigma} \cong
\left(M(\lambda)_{q^{-2i}\lambda}\right)^*$,
$\dim M(\lambda)_{q^{-2i}\lambda}^{\sigma}=1$.
For each $i \ge 0$, define $v_i^* \in M(\lambda)^*$ by
$v_i^*(v_j)=\delta_{i,j}$ (the Kr\"onecker delta)
for all $j \ge 0$.
Then $v_i^* \in M(\lambda)_{q^{-2i}\lambda}^{\sigma}$, and so
$\{v_i^* | \ i \ge 0\}$ is a $k$-basis
of $M(\lambda)^{\sigma}$.

For $i, j \ge 0$,
$(e v_i^*)(v_j)=v_i^*\left(\sigma(e)v_j\right)=v_i^*(fv_j)=v_i^*(v_{j+1})=\delta_{i,j+1}=\delta_{i-1,j}$.
Hence $e v_i^*=v_{i-1}^*$.
Similarly, $\displaystyle{f v_i^*=[i+1] \frac{q^{-i}\lambda-q^i\lambda^{-1}}{q-q^{-1}}
v_{i+1}^*} \ (i \ge 0)$. Thus $f v_i^*=0$ iff $\lambda=\pm q^i$.
Moreover, if $\lambda=\epsilon q^n \ (n \in \mathbb{Z})$,
then $fv_i^*=\epsilon [i+1][n-i]v_{i+1}^*$.

Now, utilizing the structure of $T$ modules obtained in the previous section,
we explicitly find the action of $U_q$ on
$T(n,\epsilon)^\sigma$ where $n \in \mathbb{N}_0$.
By Proposition \ref{p:DGK}, $T(n,\epsilon)^\sigma=\oplus_{i \ge 0}
T(n,\epsilon)_{\epsilon q^{n-2i}}^\sigma$ and
\begin{equation*}
 \dim T(n,\epsilon)_{\epsilon q^{n-2i}}^\sigma=
 \begin{cases}
  1,   &\mbox{if} \ 0 \le i \le n\\
  2,   &\mbox{if} \ i \ge n+1.
 \end{cases}
\end{equation*}
Recall the $k$-basis $\{v_i | \ i \ge 0\} \cup \{z_i | \ i
\ge n+1\}$ of $T(n,\epsilon)$ from Theorem \ref{t:T}.
For each $i \ge 0$, define $v_i^* \in T(n,\epsilon)^*$ by
$v_i^*(v_j)=\delta_{i,j}$ for $j \ge 0$ and
$v_i^*(z_j)=0$ for $j \ge n+1$.
Also, for each $i \ge n+1$, define $z_i^* \in T(n,\epsilon)^*$ by
$z_i^*(v_j)=0$ for $j \ge 0$ and
$z_i^*(z_j)=\delta_{i,j}$ for $j \ge n+1$.
Evidently $v_i^*$
and $z_i^*$ belong to $T(n,\epsilon)_{\epsilon q^{n-2i}}^\sigma$
and hence, by dimension count,
$\{v_i^* | \ i \ge 0\} \cup \{z_i^* | \ i \ge n+1\}$ is a $k$-basis of
$T(n,\epsilon)^\sigma$.

Since $T(n,\epsilon)$ contains a submodule spanned by $\{v_i| \ i \ge
0\}$ isomorphic to $M(\epsilon q^n)$,
we utilize the previous calculations for
$v_i^*$ restricted to $M(\epsilon q^n)$.
For $j \ge 0$, we have
$(f v_i^*)(v_j)= \epsilon[i+1][n-i] \delta_{i+1,j}$.
For $j \ge n+1$,
$ (f v_i^*)(z_j)=v_i^*(ez_j)=
v_i^*(\epsilon[j][n-j+1]z_{j-1}+v_{j-1})=v_i^*(v_{j-1})=\delta_{i,j-1}=
\delta_{i+1,j}$.
Thus, for $i \ge 0$, $f
v_i^*=\epsilon[i+1][n-i]v_{i+1}^*+z_{i+1}^*$ which is never 0. (Here, we set
$z_j=0$ and $z_i^*=0$ for $i,j < n+1$.)
Similarly, for $i \ge n+1$,
$f z_i^*=\epsilon[i+1][n-i]z_{i+1}^*$ which is also
never 0. Moreover,
$e v_i^*=v_{i-1}^*$ for $i \ge 0$ and  $e z_i^*=z_{i-1}^*$ for $i \ge n+1$.

Therefore, $T(n,\epsilon)^\sigma$ looks like
\begin{equation*}
 \begin{array}{crcll}
   \mbox{weights}           &\hspace{1in}\qquad
&\mbox{basis}           &\mbox{vectors}               &\\
   \epsilon  q^n            &\hspace{1in}\qquad                         &
                 &\cdot \ v_0^*                &\qquad \qquad ev_0^*=0\\
                            &                                           &
                 &\negthickspace\negthickspace\negthickspace e \uparrow \downarrow f      &\\
   \epsilon q^{n-2}         &\hspace{1in}\qquad                         &
                 &\cdot \ v_1^*                &\\
                            &                                           &
                 &\uparrow \downarrow          &\\
   \vdots                   &\hspace{1in}\qquad                         &
                 &\vdots              &\\
                            &                                           &
                 &\uparrow \downarrow          &\\
   \epsilon q^{-n}          &\hspace{1in}\qquad                         &
                 &\cdot \ v_n^*                &\qquad \qquad
fv_n^*=z_{n+1}^*\\
                            &                                           &f
\swarrow             &\uparrow e              &\\
   \epsilon q^{-n-2}        &\hspace{1in}\qquad z_{n+1}^* \ \cdot       &
                 &\cdot \ v_{n+1}^*            &\qquad \qquad ez_{n+1}^*=0\\
                            &\hspace{1in}\qquad e \uparrow \downarrow f &f
\swarrow\negthickspace\downarrow    &\uparrow e             &\\
   \epsilon q^{-n-4}        &\hspace{1in}\qquad z_{n+2}^* \ \cdot       &
                                &\cdot \ v_{n+2}^*      &\\
                            &\hspace{1in}\qquad \uparrow \downarrow \ \
&\swarrow\negthickspace\downarrow      &\uparrow                &\\
   \vdots                   &\hspace{1in}\qquad \vdots                  &
                                &\vdots                 &
 \end{array}
\end{equation*}

\begin{theorem}\label{t:dt}
Verma modules $M(\lambda)$ where $\lambda \ne \epsilon q^n$
for $n \in \mathbb{N}_0$ and $T$ modules $T(n,\epsilon)$
for $n \in \mathbb{N}_0$ are self-dual.
\end{theorem}

\begin{proof}
The first part of the statement follows from Proposition \ref{p:DGK} (vii).
However, it can also be seen directly
from the $U_q$-structure of $M(\lambda)^\sigma$ obtained above
since in this case
$M(\lambda)^\sigma$ is a highest weight module generated by $v_0^*$ of
weight $\lambda$ that is $U_q^-$-torsion free, and therefore
$M(\lambda)^\sigma \cong M(\lambda)$.

For the second part of the statement,
define $v_i^\di=f^iv_0^*$ for $i \ge 0$ and
$z_i^\di=\epsilon^n([n]!)^2f^{i-n-1}v_{n+1}^*$ for $i \ge n+1$. We see
from the preceeding calculations and commutation relations (\ref{commef})
that the $U_q$-action on the basis $\{v_i^\di\}_{i \ge 0} \cup \{z_i^\di\}_{i \ge n+1}$
of $T(n,\epsilon)^\sigma$
is the same as on the basis $\{v_i\}_{i \ge 0} \cup \{z_i\}_{i \ge n+1}$
of $T(n,\epsilon)$ from Theorem \ref{t:T}, and so the theorem follows.
\end{proof}

\begin{rem}
Inductively,
\begin{equation}
 v_i^\di=
 \begin{cases}
  \epsilon^i\left([i]!\right)^2 \q{n}{i} v_i^*,                  &0 \le i \le n\\
  \epsilon^{i-1}(-1)^{i-n-1}([n]![i-n-1]!)^2 \q{i}{n+1} z_i^*,   &i \ge n+1.
 \end{cases}
\end{equation}
\end{rem}

\begin{corollary} \label{Ms}
For $n \in \mathbb{N}_0$, there is a $U_q$-module isomorphism
\[
 M(\epsilon q^n)^{\sigma} \cong T(n,\epsilon) \slash M(\epsilon q^{-n-2}).
\]
\end{corollary}

\begin{proof}
By Corollary \ref{MTM}, Proposition \ref{p:DGK} (vi), and Theorem \ref{t:dt},
the following diagram of $U_q$-modules and $U_q$-module homomorphisms
\[
\begin{array}{ccccccccc}
0 & \to & M(\epsilon q^{-n-2}) & \to & T(n,\epsilon) & \to & T(n,\epsilon)/M(\epsilon q^{-n-2}) &
 \to & 0\\
  &     & \theta' \downarrow \, \cong &  & \theta \downarrow \, \cong   &     & \psi \downarrow      &
     &  \\
0 & \to & M(\epsilon q^{-n-2})^\sigma & \to & T(n,\epsilon)^\sigma & \to & M(\epsilon q^n)^\sigma
 & \to & 0
\end{array}
\]
is commutative with exact rows.
Here, $\theta : T(n,\epsilon) \to T(n,\epsilon)^\sigma$ is
a $U_q$-module isomorphism and $\theta'$ is its restriction
to $M(\epsilon q^{-n-2})$. Therefore,
the induced map
$\psi : T(n,\epsilon)/M(\epsilon q^{-n-2}) \to M(\epsilon q^n)^\sigma$
is a $U_q$-module isomorphism, as well.
\end{proof}

\begin{rem}
Utilizing Proposition \ref{FG}, we can give a different proof of the second part
of Theorem~\ref{t:dt}. Indeed, the discussion prior to Theorem \ref{t:dt} shows that
$T(n,\epsilon)^\sigma$ belongs to the category $\mathcal{I}$, i.e. $T(n,\epsilon)^\sigma$
is a weight module on which $e$ acts locally nilpotently and it is easily seen that
$f$ acts injectively. Moreover, by
Proposition \ref{p:DGK} (v) and (viii),
$T(n,\epsilon)^\sigma$ is indecomposable because $T(n,\epsilon)$ is.
Now, it follows from the weight space structure of $T(n,\epsilon)^\sigma$ and
Proposition \ref{FG} (i) that
$T(n,\epsilon)^\sigma \cong T(n,\epsilon)$.
\end{rem}

The next corollary follows immediately from Proposition \ref{FG} and Proposition \ref{p:DGK} (viii).

\begin{corollary} \label{both}
The restricted duals of modules belonging both to the category
$\mathcal{O}\left(U_q(sl_2)\right)$ and the category
$\mathcal{I}\left(U_q(sl_2)\right)$ are direct sums of modules from the set
$\{T(n,\epsilon), \, M(\epsilon q^{-n-2}), \, T(n,\epsilon)/M(\epsilon q^{-n-2}) | \, n \in \mathbb{N}_0,
\epsilon = \pm 1\}$.
\end{corollary}

\vspace{.2in}

\noindent ACKNOWLEDGMENTS \ This work is a part of the author's Ph.D. thesis
and she expresses her sincere gratitude
to her thesis advisor
Professor Vinay Deodhar for his support and valuable discussions
throughout her time as a graduate student and also for his
continued interest and constant encouragement.
Moreover, the author thanks Professor Vyjayanthi Chari for a question and the referee
for suggestions that have both led to the present version of the paper
and are greatly appreciated.

\vspace{.5in}


\end{document}